\documentclass{amsart}

\usepackage{amssymb}

\newtheorem{thm}{Theorem}

\newtheorem{cor}[thm]{Corollary}

\newtheorem{prop}[thm]{Proposition}

\theoremstyle{definition}
\newtheorem{defn}[thm]{Definition}

\newtheorem{say}[thm]{}

\newtheorem{prob}[thm]{Problem}

\newtheorem{ack}{Acknowledgments}

\newtheorem{defn-thm}[thm]{Definition--Theorem}  
\newtheorem{defn-lem}[thm]{Definition--Lemma}  

\theoremstyle{remark}

\renewcommand{\o}[0]{{\mathcal O}} 

\renewcommand{\r}[0]{{\mathbb R}} 

\renewcommand{\a}[0]{{\mathbb A}}

\newcommand{\p}[0]{{\mathbb P}}

\newcommand{\q}[0]{{\mathbb Q}}

\newcommand{\qtq}[1]{\quad\mbox{#1}\quad}

\newcommand{\mult}[0]{\operatorname{mult}}

\newcommand{\supp}[0]{\operatorname{Supp}}

\newcommand{\sing}[0]{\operatorname{Sing}}    
\newcommand{\ex}[0]{\operatorname{Ex}}

\newcommand{\cosupp}[0]{\operatorname{cosupp}}

\def\into{\DOTSB\lhook\joinrel\to}

\begin{document}
\bibliographystyle{amsalpha}

\title{Semi log resolutions}
\author{J\'anos Koll\'ar}


\maketitle

The aim of this note is to discuss resolution theorems
that are useful in the study of semi log canonical
varieties.

\begin{defn}[Simple normal crossing]\label{snc-nc.defn}
Let $k$ be a field, $X$ a  $k$-scheme and $D=\sum a_iD_i$
a Weil divisor on $X$ with the $D_i$ irreducible.

We say that $(X,D)$ has {\it simple normal crossing}
or {\it snc} at a point $p\in X$ if $X$ is smooth at $p$ 
and there are
 local coordinates  $x_1,\dots, x_{n}$
such that $\supp D\subset (x_1\cdots x_n=0)$ near $p$.
Alternatively, if  for each $D_i$ there is a
$c(i)$ such that $D_i=(x_{c(i)}=0)$ near $p$.
 
We say that $(X,D)$ has {\it  normal crossing}
or {\it nc} at a point $p\in X$
if $(\hat X_{K}, D|_{\hat X_{K}})$ is snc at $p$ where
$\hat X_{K}$ denotes the completion  at $p$
and $K$ is an algebraic closure of $k(p)$.

Let $p\in D$ be a nc point of multiplicity 2.
If the characteristic is different from 2, 
then, in suitable local coordinates, $D$
can be given by an equation $x_1^2-ux_2^2=0$ where 
$u\in \o_{p,X}$ is a unit.
$D$ is snc at $p$ iff $u$ is a square.

For example, $(y^2=x^2+x^3)\subset \a^2$
is nc but it is not snc at the origin.
Similarly, $(x^2+y^2=0)\subset \a^2$ is nc
but it is snc only if $\sqrt{-1}$ is in the base field $k$.

We say that $(X,D)$ is snc (resp.\ nc) if it is 
 snc (resp.\ nc) for every $p\in X$.

Given $(X,D)$, there is a largest open set
$U\subset X$ such that
$(U, D|_U)$ is snc (resp.\ nc).
This open set is called that {\it snc} (resp.\ {\it nc}) {\it locus}
 of $(X,D)$.

\end{defn}

\begin{defn}[Log resolution]\label{log-res.defn}
Let $k$ be a perfect field, $X$ a reduced $k$-scheme and $D$
a Weil divisor on $X$. A {\it log resolution} of
$(X,D)$ is a proper birational morphism
$f:X'\to X$ such that 
 $D':=\supp\bigl(f^{-1}(D)+\ex(f)\bigr)$
is a snc divisor on $X$. (In particular, all of its irreducible
components have codimension 1.) 
 Here
$\ex(f)$ denotes the exceptional set of $f$, that is, the set of points where
$f$ is not a local isomorphism.
We also say that  $f:(X',D')\to (X,D)$ is a log resolution.
\end{defn}

The basic existence result on resolutions was established
by \cite{hir}. We also need a strengthening of it,
due to \cite{szabo}.

\begin{thm}[Existence of log resolutions]\label{log-res.thm}
Let $X$ be an algebraic 
space of finite type over a field of characteristic 0 and $D$
a Weil divisor on $X$.
\begin{enumerate}
\item \cite{hir} $(X,D)$ has a log resolution.
\item \cite{szabo} $(X,D)$ has a log resolution $f:X'\to X$ such that
$f$ is an isomorphism over the snc locus of $(X,D)$.
\end{enumerate}
\end{thm}

 Here we show how (\ref{log-res.thm}.2)
can  be reduced to the Hironaka-type resolution theorems presented in
\cite{k-res}. The complication is that the Hironaka method and its
variants  proceed by induction on the multiplicity.
Thus, for instance, 
 the method would normally blow up every triple point of $D$ before
dealing with the non-snc double points. In the present situation, however,
we want to keep the snc triple points untouched.

 We can start by resolving the singularities of $X$,
thus it is no restriction to assume from the beginning that $X$ is smooth.
To facilitate induction, we work with a  more general
resolution problem.

\begin{defn}\label{snc.with.diff.comps}
Consider the object $(X,I_1,\dots,I_m,E)$ where $X$ is a smooth variety, 
the $I_j$ are ideal sheaves of Cartier divisors and
$E$ a snc divisor.
We say that $(X,I_1,\dots,I_m,E)$ has {\it simple normal crossing}
or {\it snc} at a point $p\in X$ if $X$ is smooth at $p$ 
and there are
 local coordinates  $x_1,\dots, x_r, x_{r+1},\dots, x_n$ and an 
injection $\sigma:\{1,\dots,r\}\to \{1,\dots,m\}$
such that 
\begin{enumerate}
\item $I_{\sigma(i)}=(x_i)$ near $p$
for  $1\leq i\leq r$ and $p\notin \cosupp I_j$ for every other $I_j$; 
\item $\supp E\subset (\prod_{i>r}x_i=0)$ near $p$.
\end{enumerate}
Thus $E+\sum_j\cosupp I_j$ has snc support near $p$,
but we also assume that no two of 
 $E, \cosupp I_1,\dots, \cosupp I_m$  have a common
irreducible component near $p$.
Furthermore, the $I_j$ are assumed to vanish with multiplicity 1,
but we do not care about the multiplicities in $E$.
The definition is chosen mainly to satisfy the
following restriction property:

\begin{enumerate}\setcounter{enumi}{2}
\item Assume that $I_1$ is the ideal sheaf
of a smooth divisor $S\subset X$ and that none of the irreducible components
of $S$ is contained in $E$ or in $\cosupp I_j$ for $j>1$. Then
$(X,I_1,\dots,I_m,E)$ is snc near $S$ iff
$(S,I_2|_S,\dots,I_m|_S,E|_S)$ is snc.\qed
\end{enumerate}

The set of all points where $(X,I_1,\dots,I_m,E)$
is  snc is  open. It is denoted by
$\operatorname{snc}(X,I_1,\dots,I_m,E)$.
\end{defn}

\begin{defn}\label{bu.seq.defn}
Let $Z\subset X$ be a smooth, irreducible subvariety
that has simple normal crossing with $E$ (cf.\ \cite[3.25]{k-res}).
Let $\pi:B_ZX\to X$ denote the  blow-up with exceptional divisor
$F\subset B_ZX$. Define
the {\it birational transform} of $(X,I_1,\dots,I_m,E)$
as
$$
(X':=B_ZX,I'_1,\dots,I'_m,E':=\pi^{-1}_{\rm tot}E)
\eqno{(\ref{bu.seq.defn}.1)}
$$
where $I'_j=g^*I_j(-F)$ if $Z\subset \cosupp I_j$
and $I'_j=g^*I_j$ if  $Z\not\subset \cosupp I_j$.
By an elementary computation, the birational
transform commutes with restriction to a smooth
subvariety (cf.\ \cite[3.62]{k-res}).
As in \cite[3.29]{k-res} we can define  blow-up sequences.
\end{defn}

The assertion (\ref{log-res.thm}.2) will be a special case of the following
 result.

\begin{prop}\label{log-res.thm.4}  
 Let $X$ be a smooth variety, $E$ an snc divisor on $X$
and $I_j$  ideal sheaves of Cartier divisors.
Then there is a smooth  blow-up sequence
$$
\Pi:(X_r, I^{(r)}_1,\dots,I^{(r)}_m,E^{(r)})\to \cdots \to 
(X_1, I^{(1)}_1,\dots,I^{(1)}_m,E^{(1)})=(X,I_1,\dots,I_m,E)
$$
such that 
\begin{enumerate}
\item $(X_r, I^{(r)}_1,\dots,I^{(r)}_m,E^{(r)})$
has snc everywhere, 
\item for every $j$, $\cosupp I^{(r)}_j$ is the
birational transform of (the closure of)
\newline $\cosupp I_j\cap \operatorname{snc}(X, I_1,\dots, I_m,E)$,   and
\item $\Pi$ is an isomorphism over 
$\operatorname{snc}(X,I_1,\dots,I_m,E)$.
\end{enumerate}
\end{prop}

Proof. The proof is by induction on
$\dim X$ and on $m$.
\medskip

{\it Step \ref{log-res.thm.4}.i.} Reduction to the case where
$I_1$ is the ideal sheaf of a smooth divisor.
\medskip

 Apply
order reduction \cite[3.107]{k-res}
to $I_1$. 
(Technically, to the marked ideal $(I_1,2)$;
see \cite[Sec.3.5]{k-res}.)
In this process, we only blow up a center $Z$ if the
(birational transform of) $I_1$ has order $\geq 2$ along $Z$.
These are contained in the non-snc locus. A slight problem is
that in  \cite[3.107]{k-res} the transformation rule used is
$I_1\mapsto \pi^*I_1(-2F)$ instead of 
$I_1\mapsto \pi^*I_1(-F)$ as in (\ref{bu.seq.defn}.1). Thus each blow-up
for  $(I_1,2)$ corresponds to two blow ups in
the sequence for $\Pi$: first we blow up
$Z\subset X$ and then we blow up $F\subset B_ZX$.

At the end the  maximal order of $I^{(r)}_1$ becomes 1.
Since $I^{(r)}_1$ is the ideal sheaf of a Cartier divisor,
 $\cosupp I^{(r)}_1$ is a disjoint union of smooth
divisors. 

\medskip

{\it Step \ref{log-res.thm.4}.ii.} Reduction to the case when
$(X,I_1, E)$ is snc.
\medskip

The first part is an easier version of Step (\ref{log-res.thm.4}.iii),
and should be read after it.
Let $S$ be an irreducible component of $E$. Write $E=S+E'$ 
 and consider the restriction
  $(S, I_1|_S,E'|_S)$.
By induction on the dimension, there is a  blow-up sequence
$\Pi_S:S_r\to\cdots\to S_1=S$ such that
$\bigl(S_r, (I_1|_S)^{(r)},(E'|_S)^{(r)}\bigr)$ is snc 
and $\Pi_S$ is an isomorphism over
$\operatorname{snc}(S, I_1|_S,E'|_S)$.
The ``same'' blow-ups give a   blow-up sequence
$\Pi:X_r\to\cdots\to X_1=X$ such that
$\bigl(X_r,I^{(r)}_1,E^{(r)}\bigr)$ is snc near $S_r$
and $\Pi$ is an isomorphism over
$\operatorname{snc}(X, I_1,E)$.

We can repeat the procedure for any other
 irreducible component of $E$. Note that as we blow up,
the new exceptional divisors are added to $E$, thus $E^{(s)}$
has more and more  irreducible components as $s$ increases.
However,  we only add new irreducible components to $E$ that
are exceptional divisors obtained by blowing up 
a smooth center that is contained in 
(the birational transform of) $\cosupp I_1$. Thus these
automatically have snc with $ I_1$. 
Therefore the procedure needs to be repeated
only for the original irreducible components of $E$.

After finitely many steps, $(X,I_1, E)$ is snc near $E$
and $X$ and $\cosupp I_1$ are smooth. Thus
$(X,I_1, E)$ is snc everywhere.

(If we want to resolve just one 
$(X, I_j,E)$, we can do these steps in any order,
but for a  functorial resolution one needs an ordering of the
index set of $E$ and proceed systematically.) 

\medskip

{\it Step \ref{log-res.thm.4}.iii.} Reduction to the case when
$(X,I_1,\dots, I_m, E)$ is snc near $\cosupp I_1$.
\medskip

Assume that $(X,I_1, E)$ is snc.
 Set $S:=\cosupp(I_1)$.
If  an irreducible component
 $S_i\subset S$ is contained in $\cosupp I_j$ for some $j>1$ then
we  blow up $S_i$. This reduces  $\mult_{S_i} I_1$ and
$\mult_{S_i} I_j$ by 1. Thus eventually none of the 
 irreducible components of
 $ S$ are contained in $\cosupp I_j$ for $j>1$.
Thus we may assume that the $I_j|_S$ are ideal sheaves of
Cartier divisors for $j>1$ and 
consider the restriction
  $(S, I_2|_S,\dots, I_m|_S,E|_S)$.

By induction there is a  blow-up sequence
$\Pi_S:S_r\to\cdots\to S_1=S$ such that
$$
\bigl(S_r, (I_2|_S)^{(r)},\dots, (I_m|_S)^{(r)},(E|_S)^{(r)}\bigr)
\qtq{is snc}
$$
and $\Pi_S$ is an isomorphism over
$\operatorname{snc}(S, I_2|_S,\dots, I_m|_S,E|_S)$.
The ``same'' blow-ups give a   blow-up sequence
$\Pi:X_r\to\cdots\to X_1=X$ such that
the restriction
$$
\bigl(S_r,I^{(r)}_2|_{S_r},\dots,I^{(r)}_m|_{S_r},E^{(r)}|_{S_r}\bigr)
\qtq{is snc}
$$
and $\Pi$ is an isomorphism over
$\operatorname{snc}(X, I_1,\dots, I_m,E)$.
(Since we use only order 1 blow-ups, this is obvious.
For higher orders, one would need 
the Going-up theorem  \cite[3.84]{k-res}, which holds only for
$D$-balanced ideals. 
Every  ideal of order 1 is $D$-balanced \cite[3.83]{k-res},
that is why we do not need to worry about subtleties here.)

As noted in (\ref{snc.with.diff.comps}.3), this implies that
$$
\bigl(X_r,\o_{X_r}(-S_r), I^{(r)}_2,\dots,I^{(r)}_m,E^{(r)}\bigr)
\qtq{is snc near $S_r$.}
$$
Note, furthermore, that $S_r=\cosupp I^{(r)}_1$, 
hence 
$\bigl(X_r,I^{(r)}_1,\dots,I^{(r)}_m,E^{(r)}\bigr)$ is snc near 
$\cosupp I^{(r)}_1$.

\medskip

{\it Step \ref{log-res.thm.4}.iv.} Induction on $m$.
\medskip

By Step 3, we can assume that $(X,I_1,\dots, I_m, E)$ is snc 
near $\cosupp I_1$.
Apply (\ref{log-res.thm.4})  to
$(X,I_2,\dots, I_m, E)$. 
The resulting $\Pi:X_r\to X$ is an isomorphism
over $\operatorname{snc}(X, I_2,\dots, I_m,E)$.
Since $\cosupp I_1$ is contained in
 $\operatorname{snc}(X, I_2,\dots, I_m,E)$,  all
the blow up centers are disjoint from $\cosupp I_1$.
Thus $\bigl(X_r,I^{(r)}_1,\dots,I^{(r)}_m,E^{(r)}\bigr)$
is also snc.

Finally, we may blow up any irreducible component of
$\cosupp I^{(r)}_j$ that is not the birational transform of
an irreducible component of $\cosupp I_j$
which intersects  $\operatorname{snc}(X, I_1,\dots, I_m,E)$.
\qed

\begin{say}[Proof of (\ref{log-res.thm})]
Let $D_j$ be the irreducible components of
$D$. Set $I_j:=\o_X(-D_j)$ and $E:=\emptyset$.
Note that $(X,D)$ is snc at $p\in X$ iff
$(X, I_1,\dots, I_m,E)$  is snc at $p\in X$.

If $X$ is a variety, we can apply (\ref{log-res.thm.4}) to
$(X, I_1,\dots, I_m,E)$ to get
$\Pi:X_r\to X$ and
$(X_r, I^{(r)}_1,\dots,I^{(r)}_m,E^{(r)})$. 
Note that $E^{(r)}$ contains the whole exceptional
set of $\Pi$, thus the support of
$D'=\Pi^{-1}_*D+\ex(\Pi)$
is contained in
$E^{(r)}+\sum_j \cosupp I^{(r)}_j$. Thus
$D'$ is snc. By (\ref{log-res.thm.4}.3),
$\Pi$ is an isomorphism over
the snc locus of $(X,D)$.

The resolution constructed in (\ref{log-res.thm.4})
commutes with smooth morphisms and with change of fields
\cite[3.34.1--2]{k-res}, at least if
in (\ref{bu.seq.defn}) we allow   reducible blow-up centers.

 As in \cite[3.42--45]{k-res}, we conclude that
(\ref{log-res.thm}) and (\ref{log-res.thm.4})
also hold for algebraic and analytic spaces over
a field of characteristic 0.

Starting with $(X,D)$, the above proof depends on 
an ordering of the irreducible components of $D$.
This is an artificial device, but I don't know how to
avoid it. \qed
\end{say}

\begin{say}\label{pinch.problem.say} It should be noted that
(\ref{log-res.thm}.2) fails for nc instead of snc.
The simplest example is given by the {\it pinch point}
$D:=(x^2=y^2z)\subset \a^3=:X$.
Here $(X,D)$ has nc outside the origin.
At a point along
the $z$-axis, save at the origin,
$D$ has 2 local analytic branches. As we go around the origin,
these 2 branches are interchanged.
This continues to hold after any birational map that
is an isomorphism over the generic point of the $z$-axis
and so we can never get rid of the pinch point
without blowing up the $z$-axis.

Note that $(X,D)$  is not snc along the $z$-axis,
thus in constructing a log resolution as in (\ref{log-res.thm}.2),
we are allowed to blow up the $z$-axis. 

This leads to the following general problem:

\begin{prob} Describe the smallest class of singularities
${\mathcal  S}$ such that for every 
$(X,D)$ there is a proper birational map $f:X'\to X$
such that 
\begin{enumerate}
\item $(X',D')$ has only singularities in ${\mathcal  S}$, and
\item $f$ is an isomorphism over the nc locus of $(X,D)$.
\end{enumerate}
\end{prob}

In dimension 2 we can take, up to \'etale equivalence, 
${\mathcal  S}= \{(xy=0)\subset \a^2\}$
and in dimension 3 we can almost certainly take
$$
{\mathcal  S}= \{(xy=0), (xyz=0), (x^2=y^2z)\subset \a^3\}.
$$

In higher dimensions, there is not even a clear conjecture on what
${\mathcal  S}$ should be.
\end{say}

\begin{defn}[Semi snc] The ideal local model of an snc
 $\q$-divisor is given by
$D=\sum_{i=1}^n a_i (x_i=0)$ on $X=\a^n$. 
We can also view this as sitting on $\a^{n+1}$,
where $X=(x_{n+1}=0)$ and $D$ is defined using the other coordinates.

Following this example, we can define a non-normal version of
snc where $X\subset \a^{n+1}$ is defined by the product of some
of the coordinates and $D$ is defined using the remaining coordinates.

For  $n=2$  we get three possible local models.
\begin{enumerate}
\item $S=(z=0)\subset \a^3$ and $D=a_x(x|_S=0)+a_y(y|_S=0)$.
This is the usual normal case.
\item $S=(yz=0)\subset \a^3$ and  $D=a_x(x|_S=0)$. 
Note that as a Weil divisor, $D$ has two irreducible
components, namely $D_1:=(x=y=0)$ and $D_2:=(x=z=0)$. 
The  support of the Weil $\r$-divisor $a_1D_1+a_2D_2$
is always snc, but the pair
$(S, a_1D_1+a_2D_2)$ is semi-snc only if $a_1=a_2$. 
It is easy to see that $a_1D_1+a_2D_2$ is $\r$-Cartier only if
$a_1=a_2$. 
\item $S=(xyz=0)\subset \a^3$ and  $D=0$.
\end{enumerate}

Based on this, local models of {\it semi-snc}
pairs are the following.

Let $Y$ be a smooth variety, $0\in Y$ a point
and $y_1,\dots, y_{n+1}$ local coordinates.
Let $I_X, I_D\subset \{1,\dots,n+1\}$
be disjoint subsets and $c:I_D\to \r$ a function.

Then 
$$
X:=\sum_{i\in I_X} (y_i=0)=\sum_{i\in I_X} X_i
$$
is an snc divisor on $Y$, which we view now as a subscheme, and
$$
D:=\sum_{i\in I_D} c(i) \bigl( y_i|_X=0\bigr)=\sum_{i\in I_D} c(i) D_i
$$
is a Weil  $\r$-divisor on $X$.

Let $X$ be a reduced variety and
$D$ a Weil $\q$-divisor on $X$. We say that
$(X,D)$ is {\it semi-snc}
if every point $x\in X$ has an open neighborhood
$x\in U$ such that $(U,D|_U)$ is isomorphic to
a local model constructed above.

As in (\ref{snc-nc.defn}), one can also define
{\it semi-nc}.
\end{defn}

\begin{say}[Semi log resolutions]
What is the right notion of resolution or log resolution for
non-normal varieties?

The simplest choice is to make no changes and work with resolutions.
In particular, if $X=\cup_i X_i$ is a reducible scheme
and $f:X'\to X$ is a resolution then $X'=\cup_i X'_i$
such that each $X'_i\to X_i$ is a resolution.
Note that we have not completely forgotten the gluing
data determining $X$ since $f^{-1}(X_i\cap X_j)$ is
part of the exceptional set, and so we keep track of it.

There are, however, several inconvenient aspects.
For instance, $f_*\o_{X'}\neq \o_X$, and this makes it
difficult to study the Picard group of $X$ or the cohomology of
line bundles on $X$ using $X'$. Another problem is that
although $\ex(f)$ tells us which part of $X_i$ intersects
the other components, it does not tell us anything about
what the actual isomorphism is between
$(X_i\cap X_j)\subset X_i$ and $(X_i\cap X_j)\subset X_j$.

It is not clear how to remedy these problems for an arbitrary
reducible scheme, but we are dealing with with schemes
that have only double normal crossing in codimension 1.

We can thus look for $f:X'\to X$ such that
$X'$ has only double normal crossing singularities and
 $f$ is an isomorphism over codimension 1  points of $X$.

As in (\ref{log-res.thm}), this works for simple nc
but not in general. We need to allow at least pinch points.
\end{say}

\begin{defn}[Pinch points]\label{pinch.pt.defn}
 Let $X$ be a smooth variety
over a field of characteristic $\neq 2$ and $D\subset X$ a divisor.
We say that $D$ has a {\it pinch point} at $p\in D$
if, in suitable local coordinates, $D$ can be defined
by the equation $x_1^2-x_2^2x_3=0$. 

Note that this notion is invariant under field extensions
and even completion. Indeed,  if the singular set of
$D$ is a codimension 2 smooth subvariety, then
$D$ can be locally given by an equation
$a x_1^2+bx_1x_2+cx_2^2=0$ where $a,b,c$ are regular functions.
If the quadratic part of the equation is a  square times a unit,
then, after a coordinate change,  we can write the equation as
 $x_1^2+cx_2^2=0$.  This gives a pinch point after a field extension
and  completion iff the linear term of $c$ is
independent of $x_1, x_2$. Thus we can take
$x_3=-c$ to get the  equation $x_1^2-x_2^2x_3=0$.

Let us blow up $Z:=(x_1=x_2=0)$. The normalization of $D$ is contained
in  the affine charts with coordinates
$x'_1:=x_1/x_2, x_2,\dots, x_n$.
If we introduce $x'_3:=x_3-{x'_1}^2$ then the normalization of $D$ is
given by $(x'_3=0)$. The preimage of $Z$ is the smooth divisor
$x_2=0$ and the involution on it is 
$(x'_1, 0,0, x_4,\dots, x_n)\mapsto (-x'_1, 0,0, x_4,\dots, x_n)$.

A function $f$ defines a $\tau$-invariant divisor
iff
$$
f(x'_1,x_2, x_4,\dots, x_{n})=
\left\{
\begin{array}{l}
\hphantom{x'_1}
g({x'_1}^2, x_4,\dots, x_{n})+x_2h(x'_1, x_2,x_4,\dots, x_{n}),\qtq{or}\\
x'_1g({x'_1}^2, x_4,\dots, x_{n})+x_2h(x'_1, x_2,x_4,\dots, x_{n}).
\end{array}
\right.
$$
In the first case $f$ is $\tau$-invariant and descends to a
regular function on $D$.
In the second case $f$ is not $\tau$-invariant, but $f^2$  descends to a
regular function on $D$.

In particular, 
$(x_1=x_3=0)\subset (x_1^2=x_2^2x_3)$ is not a Cartier divisor but
 it is $\q$-Cartier since
$2(x_1=x_3=0)=(x_3=0)$ is Cartier.
\end{defn}

\begin{thm} \label{sres.nc2.thm}
Let $X$ be a  reduced scheme over a field of characteristic 0.
 Let $X^{ncp}\subset X$ be an open subset such that
$X^{ncp}$ has only smooth points  $(x_1=0)$,
double nc points  $(x_1^2-ux_2^2=0)$ and pinch points
$(x_1^2=x_2^2x_3)$.
Then there is a projective birational morphism
 $f:X'\to X$ such that
\begin{enumerate}
\item $X'$ has only smooth points,
double nc points  and pinch points,
\item $f$ is an isomorphism over $X^{ncp}$,
\item $\sing X'$ maps birationally onto the
closure of $\sing X^{ncp}$.
\end{enumerate}
\end{thm}

If $X'$ has any pinch points then they are
on an irreducible component of $B\subset \sing X'$
along which $X'$ is nc but not snc.
Then,  by (\ref{sres.nc2.thm}.3), $X$ is nc but not snc along
$f(B)$. Thus we obtain the following simple nc version.

\begin{cor}\label{sres.nc2.cor}
 Let $X$ be a  reduced scheme over a field of characteristic 0.
 Let $X^{snc2}\subset X$ be an open subset which
 has only smooth points  $(x_1=0)$ and 
simple  nc points  of multiplicity $\leq 2$ $(x_1x_2=0)$.
Then there is a projective birational morphism
 $f:X'\to X$ such that
\begin{enumerate}
\item $X'$ has only smooth points and simple 
 nc points   of multiplicity $\leq 2$,
\item $f$ is an isomorphism over $X^{snc2}$,
\item $\sing X'$ maps birationally onto the
closure of $\sing X^{snc2}$. \qed
\end{enumerate}
\end{cor}

\begin{say}[Proof of (\ref{sres.nc2.thm})] \label{pf.of.sres.nc2.thm}
The method of \cite{hir} reduces the multiplicity
of a scheme starting with the highest multiplicity locus.
We can use it to find a proper birational morphism
 $g_1:X_1\to X$ such that every point of
$X_1$ has multiplicity $\leq 2$ and $g_2$ is an
isomorphism over $X^{ncp}$. Thus by replacing $X$ by $X_1$ we may assume
to start with that every point of
$X$ has multiplicity $\leq 2$.

The next steps of the Hironaka method would not distinguish the
nc locus (that we want to keep intact) from the
other multiplicity 2 points (that we want to eliminate).
Thus we proceed somewhat differently.

Let $n:\bar X\to X$ be the normalization with reduced conductor
$\bar B\subset \bar X$.

Near any point of $X$, in local analytic or \'etale coordinates
we can write $X$ as
$$
X=\bigl(y^2=g({\mathbf x})h({\mathbf x})^2\bigr)\subset \a^{n+1}
$$
where $({\mathbf x}):=(x_1,\dots,x_n)$
and $g$ has no multiple factors.
(We allow $g$ and $h$ to have common factors.)
The normalization is then given by
$$
\bar X=\bigl(z^2=g({\mathbf x})\bigr)\qtq{where $z=y/h({\mathbf x})$.}
$$
Here $\bar B=(h({\mathbf x})=0)$ and
the involution $\tau:(z,{\mathbf x})\mapsto (-z, {\mathbf x})$
is well defined on $\bar B$.
(By contrast, the $\tau$ action on $\bar X$ depends on the
choice of the local coordinate system.)

Thus we have a pair $(Y_2,B_2):=(\bar X, \bar B)$
plus an involution $\tau_2:B_2\to  B_2$ such that
for every $b\in  B_2$ there is an \'etale neighborhood
$U_b$ of $\{b,\tau_2(b)\}$ such that
$\tau_2$ extends (nonuniquely) to an involution
$\tau_{2b}$ of $(U_b, {B_2}|_{U_b})$.

Let us apply an \'etale local resolution procedure
(as in \cite{wlod} or
\cite{k-res})
to $(Y_2,B_2)$. Let the first blow up center be
 $Z_2\subset Y_2$.
Since the procedure is \'etale local, we see that
$U_b\cap Z_2$ is $\tau_{2b}$-invariant for every $b\in B_2$.
Let $Y_3\to Y_2$ be the blow up of $Z_2$ and let
$B_3\subset Y_3$ be the birational transform of $B_2$.
Then $\tau_2$ lifts to an involution
$\tau_3$ of $B_3$ and the $\tau_{2b}$ lift to extensions 
on suitable neighborhoods.
Moreover, the exceptional divisor of $Y_3\to Y_2$ intersected
with $B_3$ is
$\tau_3$-invariant. In particular, there is an ample
line bundle $L_3$ on $Y_3$ such that $L_3|_{B_3}$ is
$\tau_3$-invariant.

At the end we obtain $g:Y_r\to Y_2=\bar X$ such that
\begin{enumerate}
 \item $Y_r$ is smooth and  $\ex(g)+B_r$ is an snc divisor,
\item $B_r$ is smooth and $\tau$ lifts to an involution $\tau_r$ on $B_r$, and
\item there is a $g$-ample line bundle $L$ such that
$L|_{B_r}$ is $\tau_r$-invariant.
\end{enumerate}
\noindent The fixed point set of $\tau_r$ is a disjoint union of smooth
subvarieties of $B_r$. By blowing up those components whose dimension
is $<\dim B_r-1$, we also achieve (after replacing $r+1$ by $r$) that
\begin{enumerate}\setcounter{enumi}{3}
 \item the fixed point set of $\tau_r$ has pure codimension 1 in $B_r$.
\end{enumerate}

Let $Z_r:=B_r/\tau_r$ and $X_r$ the universal
push-out of $Z_r\leftarrow B_r\into Y_r$, cf.\ \cite[Thm.3.1]{artin70}.

Further, let $D$ be a divisor on $Y_r$ such that
$D|_{B_r}$ is $\tau_r$-invariant. As noted in (\ref{pinch.pt.defn}),
$2D$ is the pull back of a Cartier divisor on $X_r$.
In particular, if $D$ is ample then $X_r$ is
projective. \qed
\end{say}
\medskip

We would like not just a semi resolution of $X$ but a log resolution
of the pair $(X,D)$. Thus we need to take into account the
singularities of $D$ as well. As we noted in (\ref{pinch.problem.say}),
this is not obvious even when $X$ is a smooth 3-fold.
The following weaker version, which gives the expected result
 only for the codimension 1
part of the singular set of $(X,D)$, will be sufficient for us.

\begin{thm} \label{logres.nc.thm}
Let $X$ be a  reduced scheme over a field of characteristic 0 
and $D$
a Weil divisor on $X$. Let $X^{nc2}\subset X$ be an open subset which
 has only nc points of multiplicity $\leq 2$ and 
 $D|_{X^{nc2}}$ is smooth and disjoint from
$\sing X^{nc2}$.
Then there is a projective birational morphism
 $f:X'\to X$ such that
\begin{enumerate}
\item the local models for $\bigl(X', D':=f^{-1}_*(D)+\ex(f)\bigr)$ are
\begin{enumerate}
\item (Smooth)  $X'=(x_1=0)$ and $D'=(\prod_{i\in I} x_i=0)$
for some $I\subset \{2,\dots,n+1\}$, 
\item (Double nc)  $X'=(x_1^2-ux_2^2=0)$ and $D'=(\prod_{i\in I} x_i=0)$
for some $I\subset \{3,\dots,n+1\}$, or 
\item (Pinched)  $X'=(x_1^2=x_2^2x_3)$ and $D'=(\prod_{i\in I} x_i=0)+D_2$
for some $I\subset \{4,\dots,n+1\}$
where either $D_2=0$ or $D_2=(x_1=x_3=0)$.
\end{enumerate}
\item $f$ is an isomorphism over $X^{nc2}$.
\item $\sing X'$ maps birationally onto the
closure of $\sing X^{nc2}$.
\end{enumerate}
\end{thm}

As before, (\ref{logres.nc.thm}) implies the simple nc version:

\begin{cor}\label{logres.snc.cor}
 Let $X$ be a  reduced scheme over a field of characteristic 0 
and $D$
a Weil divisor on $X$. Let $X^{snc2}\subset X$ be an open subset 
which  has only snc points of multiplicity $\leq 2$ and 
 $D|_{X^{snc2}}$ is smooth and disjoint from
$\sing X^{snc2}$.
Then there is a projective birational morphism
 $f:X'\to X$ such that
\begin{enumerate}
\item the local models for $\bigl(X', D':=f^{-1}_*(D)+\ex(f)\bigr)$ are
\begin{enumerate}
\item (Smooth)  $X'=(x_1=0)$ and $D'=(\prod_{i\in I} x_i=0)$
for some $I\subset \{2,\dots,n+1\}$, or
\item (Double snc)  $X'=(x_1x_2=0)$ and $D'=(\prod_{i\in I} x_i=0)$
for some $I\subset \{3,\dots,n+1\}$.
\end{enumerate}
\item $f$ is an isomorphism over $X^{snc2}$.
\item $\sing X'$ maps birationally onto the
closure of $\sing X^{snc2}$.
\end{enumerate}
\end{cor}

\begin{say}[Proof of (\ref{logres.nc.thm})]
 First we use (\ref{sres.nc2.thm}) to reduce to the case
when $X$ has only double nc and pinch points.
Let $\bar X\to X$ be the normalization and
$\bar B\subset \bar X$ the conductor.
Here $\bar X$ and $\bar B$ are both smooth.

Next we want to apply embedded resolution to
$(\bar X, \bar B+\bar D)$.
One has to be a little careful with $D$ since the preimage
$\bar D\subset \bar X$ need not be $\tau$-invariant.

As a first step, we 
move the support of $\bar D$ away from $\bar B$. As in
\cite[3.102]{k-res} this is equivalent to
multiplicity reduction for a suitable ideal
$I_D\subset \o_{\bar B}$. Let us now apply
multiplicity reduction for the  ideal
$I_D+ \tau^*I_D$. All the steps are now
$\tau$-invariant, so at the end we obtain
$g:Y_r\to \bar X$ such that
$B_r+D_r+\ex(g)$ has only snc along $B_r$ and
$\tau$ lifts to an involution $\tau_r$.

As in the proof of (\ref{sres.nc2.thm}), we can also assume that
the fixed locus of $\tau_r$ has pure codimension 1
in $B_r$ and that there is a $g$-ample line bundle
$L$ such that $L|_{B_r}$ is $\tau_r$-invariant.

As in the end  of (\ref{pf.of.sres.nc2.thm}),
let  $X_r$ be the universal
push-out of $B_r/\tau_r\leftarrow B_r\into Y_r$.
Then  
$(X_r, D'_r)$ has the required normal form
along  $\sing X_r$. The remaining singularities of
$D'_r$ can now be resolved as in (\ref{log-res.thm}).\qed
\end{say}

The following analog of (\ref{log-res.thm}) is still open:

\begin{prob} Let $X$ be a  reduced scheme over a field of characteristic 0 
and $D$
a Weil divisor on $X$. Let $X^{snc}\subset X$ be the largest
 open subset such that
$(X^{snc}, D|_{X^{snc}})$ is semi snc.
Is there  a projective birational morphism
 $f:X'\to X$ such that
\begin{enumerate}
\item  $(X', D')$ is semi snc and
\item $f$ is an isomorphism over $X^{snc}$?
\end{enumerate}
\end{prob}

The following weaker version is sufficient for many applications.
We do not guarantee that $f:X'\to X$ is an isomorphism over $X^{snc}$,
only that $f$ is an isomorphism over an open subset  $X^0\subset X^{snc}$
that intersects every semi log canonical center of $(X^{snc}, D|_{X^{snc}})$.
(One can see easily that the latter are exactly the
irreducible
components of  intersections of irreducible
components of $X^{snc}$ and of $D|_{X^{snc}}$.)
This implies that we do not introduce
any ``unnecessary''  $f$-exceptional divisors with
discrepancy $-1$. The latter  is usually the key property that one needs.

Unfortunately, the proof only works in the quasi projective case.

\begin{prop}\label{sdlt.weaksz.prop} Let $X$ be a  reduced quasi projective
scheme over a field of characteristic 0 
and $D$
a Weil divisor on $X$. Let $X^0\subset X$ be an open subset such that
$(X^0, D|_{X^0})$ is semi snc.
There is a projective birational morphism
 $f:X'\to X$ such that
\begin{enumerate}
\item  $(X', D')$ is semi snc and
\item $f$ is an isomorphism over the generic point of every
semi log canonical center of $(X^0, D|_{X^0})$.
\end{enumerate}
\end{prop}

Proof. In applications it frequently happens that
$X+B$ is a divisor on a variety $Y$
and $D=B|_X$. Applying (\ref{log-res.thm}.2)
to $(Y,X+B)$ gives (\ref{sdlt.weaksz.prop}).
In general, not every $(X,D)$ can be obtained this way,
but one can achieve something similar at the price of
introducing other singularities.

Take an embedding $X\subset \p^N$.
Pick a finite set $W\subset X$ such that  each
semi log canonical center of $(X^0, D|_{X^0})$
contains a point of $W$.  

Choose $d\gg 1$ such that the scheme theoretic
base locus of $\o_{\p^N}(d)(-X)$ is $X$ near every point
of $W$. Taking a complete intersection of 
$(N-\dim X-1)$ general members
in $|\o_{\p^N}(d)(-X)|$, we obtain $Y\supset X$
such that $Y$ is  smooth  at every point of $W$.
(Here we use that $X$ has only hypersurface singularities
near $W$.)

For every $D_i$ choose $d_i\gg 1$
such that the scheme theoretic
base locus of $\o_{\p^N}(d_i)(-D_i)$ is $D_i$ near every point
of $W$. For each $i$, let $D_i^Y\in |\o_{\p^N}(d_i)(-D_i)|$
be a general member.

We have thus constructed  a pair $(Y, X+\sum D_i^Y)$ 
such that 
\begin{enumerate}
\item $(Y, X+\sum D_i^Y)$ is snc near $W$, and
\item $(X, \sum D_i^Y|_X)$ is isomorphic to
$(X, \sum D_i)$ in a suitable neighborhood
of $W$.
\end{enumerate}

By (\ref{log-res.thm}.2) there is a semi log resolution of 
$$
f:(Y', X'+\textstyle{\sum} B_i)\to (Y, X+\textstyle{\sum} D_i^Y)
$$
such that $f$ is an isomorphism over an
open neighborhood of $W$.
Then $f|_{X'}:X'\to X$ is the 
log resolution we want.
\qed

 \begin{ack}
Partial financial support  was provided by  the NSF under grant number 
DMS-0758275.
\end{ack}

\bibliography{refs}

\vskip1cm

\noindent Princeton University, Princeton NJ 08544-1000

\begin{verbatim}kollar@math.princeton.edu\end{verbatim}

\end{document}